\def\nin{\noindent}

\documentclass[12pt]{article}
\usepackage{a4,amsmath,amssymb,harvard,theorem,time}
\usepackage[T1]{fontenc}
\usepackage{mathptmx}
\harvardparenthesis{square}

\date{\today\ \now}
\theoremstyle{break}
\newtheorem{Thm}{Theorem}
\theorembodyfont{\itshape}
\theoremheaderfont{\scshape}
\def\R{\mathbb R}
\newcommand{\sn}{\mathop{\textrm{sn}}\nolimits}
\newcommand{\cn}{\mathop{\textrm{cn}}\nolimits}
\newcommand{\dn}{\mathop{\textrm{dn}}\nolimits}
\def\today{\number\year\space
 \ifcase\month\or January\or February\or March\or April\or May\or June\or
   July\or August\or September\or October\or November\or December\fi
     \space\number\day
}

\title{\vspace{-15mm}\bf Another universal differential equation}

\author{Keith Briggs,
BTexact Technologies\\
Adastral Park, Antares 2 pp5\\
Suffolk IP5 3RE, UK\\
{\tt Keith.Briggs@bt.com}
}

\date{\today} %

\begin{document}
\maketitle

\begin{abstract}
\noindent
I construct a new universal differential equation (\ref{B}), in the sense 
of Rubel.  That is, its solutions approximate to arbitrary accuracy any 
continuous function on any interval of the real line.
\end{abstract}

\section{Introduction}

\cite{rubel81} proved the surprising theorem:
\begin{Thm}
Given any continuous function $\phi:\R\rightarrow\R$ and any positive
continuous function $\epsilon:\R\rightarrow\R_{+}$, there exists
a $C^\infty$ solution $y$ of
\begin{multline}
3{y'}^4y''{y''''}^2-4{y'}^4{y'''}^2{y''''}+6{y'}^3{y''}^2{y'''}{y''''}+\\
24{y'}^2{y''}^4y''''-12{y'}^3{y''}{y'''}^3-29{y'}^2{y''}^3{y'''}^2+12{y''}^7=0
\tag{R}\label{R}\nonumber
\end{multline}
such that
$
|y(t)-\phi(t)|<\epsilon(t)\qquad \forall\ t\in \R.
$
\end{Thm}
Note that Rubel's function $y$ is $C^\infty$ but not real-analytic, 
typically having a countable number of essential singularities.
The solutions of (\ref{R}) have the additional properties that they 
may be taken to be monotone if $\phi$ is monotone, and can be made to
agree with $\phi$ at a countable number of distinct abscissae.
We may rephrase the result by saying that all continuous functions
are uniform limits of sequences of solutions of (\ref{R}), or that
solutions of (\ref{R}) are dense in the space $C(\R)$ of continuous functions.

By a {\sl universal ODE\/} I will mean any $n$th order homogeneous non-trivial
polynomial ODE $P(y,y',y'',\dots,y^{(n)})=0$ (that is, an algebraic 
differential equation or ADE) whose solutions have the same property 
as the function $y$ in Rubel's theorem.    
Universal differential equations are of interest
in the theory of computability by analog computers.   For discussions
of this field, see \cite{boshernitzan85,rubel83,rubel92,rubel96}.

There have been several recent developments, which I now describe.
\cite{duffin81} found analogous universal differential equations whose
solutions are $C^n$ for finite $n>3$:
\begin{Thm}
The differential equations
\begin{equation} %
n^2{y''''}{y'}^2+3n(1-n){y'''}{y''}{y'}+(2n^2-3n+1){y''}^3=0 \tag{D1}\label{D1}\nonumber
\end{equation}
and
\begin{equation} %
n{y''''}{y'}^2+(2-3n){y'''}{y''}{y'}+2(n-1){y''}^3=0 \tag{D2}\label{D2}\nonumber
\end{equation}
are universal.
\end{Thm}
(Note that \cite{duffin81} contains several errors which I have corrected
in the above statement of (\ref{D1}).)

\cite{boshernitzan86} gave a general survey of the field and proved a new result
for finite intervals without a smoothness assumption:

\begin{Thm}
There exists an ADE (of order less than 20), the polynomial solutions
of which are dense in $C(I)$ for any finite interval $I$.
\end{Thm}
The ADE was not explicitly constructed.  Of course, by Weierstrass' classical
theorem on approximation by polynomials, we know that the set of all
polynomials $P(I)$ on an interval $I$ is dense in $C(I)$; 
the interest in theorems such as the above is that a subset of $P(I)$ is still dense.

\cite{elsner92} proved a new result for Lipschitz functions, and
\cite{elsner99} constructed a universal functional equation.
(This paper also contains an error: the statement of Theorem 1.1 should
say that the greatest coefficient is 39813120000, not 35831890000.)

\section{The new universal differential equation}

\begin{Thm}
\begin{equation}
{y''''}{y'}^2-3{y'''}{y''}{y'}+2(1-n^{-2}){y''}^3=0 \tag{B}\label{B}\nonumber
\end{equation}
is universal for $n>3$.
\end{Thm}
{\sc\footnotesize PROOF:}
I will make use of the Jacobian elliptic functions which depend on a real 
parameter $m \in \left[0,1\right)$ and satisfy
\begin{eqnarray}
\sn'(x,m)&=&\cn(x,m)\dn(x,m)\nonumber\\
\cn'(x,m)&=&-\sn(x,m)\dn(x,m)\nonumber\\
\dn'(x,m)&=&-m\sn(x,m)\cn(x,m)\nonumber\\
\sn^2(x,m)&+&\cn^2(x,m)=1\nonumber\\
\dn^2(x,m)&+&m\sn^2(x,m)=1.\nonumber
\end{eqnarray}
These properties alone will suffice for the proof.
I let $g(x)=\cn^n(x,m)$ and consider the expression
\begin{equation}
ng'''(x)g^2(x)+bg''(x)g'(x)g(x)+c{g'(x)}^3. \tag{*}\label{star}\nonumber
\end{equation}
We wish to find (if they exist) values of $m,n,b,c$ which make the expression
(\ref{star}) identically zero.  After substitution of the expression for $g(x)$,
(\ref{star}) simplifies to
\begin{eqnarray}
&&(1-m)\left[2-b+(b+c-3)n+n^2\right]+\nonumber\\
&&(2m-1)n(b+c+n)\cn^2(x,m)-\nonumber\\
&&m\left[2+b+(b+c+3)n+n^2\right]\cn^4(x,m)\nonumber  
\end{eqnarray}
and we have the solutions 
$$
\begin{tabular}{lll}
$m=0$,   & $b=2-3n$, & $c=2(n-1);   $\nonumber\\
$m=1/2$, & $b=-3n$,  & $c=2(n-1/n); $\nonumber\\
$m=1$,   & $b=-2-3n$,& $c=2(n+1)    $\nonumber.
\end{tabular}
$$
The first solution in fact reduces to (\ref{D1}); the 
last is excluded as the elliptic functions are not defined for $m=1$,
leaving only the case $m=1/2$ to consider.
The function $g$ and its first $n-2$ derivatives 
vanish when $x=\pm K(m)$, where $K$ is the complete
elliptic integral of the first kind.  Thus, if we extend $g$ to the
whole real line by making it zero outside $[-K(m),K(m)]$, the extended
function is still a solution of (\ref{B}).   Also, affine transformations
of $g$ of the
form $g(\alpha x+\beta)$ are also solutions by homogeneity, and so are
sums of such functions on disjoint intervals.
Thus
\begin{eqnarray}
Y(x)=\gamma \int_{-\infty}^{x} g(\alpha t+\beta)\text{d}t +\delta\nonumber
\end{eqnarray}
forms an S-module in the sense of \cite{rubel81} on $[K(m),K(m)]$.
(That is, it is monotonic and its first $n-1$ derivatives vanish at the 
endpoints.)  Our complete solution $y$ is then obtained by pasting together 
sufficiently many such S-modules.
$\blacksquare$

Note that (\ref{D1}), (\ref{D2}), and (\ref{B}) all reduce to the same ADE in the 
limit $n\rightarrow\infty$.
It would be interesting to know which other ADEs of the general form of (\ref{B}) are universal.
\clearpage
\nocite{*}
\bibliography{universal-ode}

\small
\vfill \nin{\tt tantalum.bt-sys.bt.co.uk:tex/universal-ode.tex} \hfill Typeset in \LaTeXe
\end{document}